\documentclass[10pt]{article}

\usepackage{graphics}
\usepackage{graphicx}

\usepackage[a4paper, left=35mm,right=35mm,top=34mm,bottom=34mm]{geometry}
\usepackage[utf8]{inputenc}
\usepackage[T1]{fontenc}
\usepackage[english]{babel}

\usepackage{enumerate}
\usepackage{graphicx}
\usepackage{hyperref}
\hypersetup{
    colorlinks=true,
    linkcolor=blue,
    filecolor=magenta,      
    urlcolor=cyan,
}
\usepackage{listings}
\usepackage{color}

\definecolor{dkgreen}{rgb}{0,0.6,0}
\definecolor{gray}{rgb}{0.5,0.5,0.5}
\definecolor{mauve}{rgb}{0.58,0,0.82}
\lstdefinelanguage{MRGC++}{%
  language=C++,
  morekeywords={T, U, MPI_Irecv, MPI_Isend, MPI_Allreduce, MPI_Waitall, Compute, Map, abs, max, Swap, MPI_Recv_init, MPI_Send_init, MPI_Startall, Copy, Init, InitRecv, InitSend, InitAllReduce, Send, Recv, AllReduce, Finalize, InitSnapshot, Snapshot, SwitchAsync, SnapReduce, MPI_Test, MPI_Start}
}
\usepackage{mathtools,amsthm,amssymb,amsfonts}
\usepackage{algorithm}
\usepackage{algorithmic}
\makeatother
\theoremstyle{plain}
\newtheorem{theorem}{Theorem}
\newtheorem{corollary}{Corollary}
\newtheorem{proposition}{Proposition}

\newtheorem{lemma}{Lemma}
\theoremstyle{definition}
\newtheorem{definition}[theorem]{Definition}
\theoremstyle{definition}

\theoremstyle{remark}
\newtheorem*{remark}{Remark}

\usepackage{caption} 
\usepackage{fancyhdr}

\author{
  {\normalsize Guillaume Gbikpi-Benissan}\thanks{Universit\'e Paris-Saclay, CentraleSup\'elec, Gif-sur-Yvette, France
    (guibenissan@gmail.com).}
  \and
  {\normalsize Qinmeng Zou}\thanks{Universit\'e Paris-Saclay, CentraleSup\'elec, Gif-sur-Yvette, France
    (zouqinmeng@hotmail.com).}
  \and
  {\normalsize Fr\'ed\'eric Magoul\`es}\thanks{Universit\'e Paris-Saclay, CentraleSup\'elec, Gif-sur-Yvette, France
    (correspondence, frederic.magoules@hotmail.com).}
}
\title{Asynchronous iterations of HSS method for non-Hermitian linear systems}
\date{}

\begin{document}
\maketitle
\thispagestyle{fancy}

\begin{abstract}
\noindent A general asynchronous alternating iterative model is designed, for which convergence is theoretically ensured both under classical spectral radius bound and, then, for a classical class of matrix splittings for $\mathsf H$-matrices. The computational model can be thought of as a two-stage alternating iterative method, which well suits to the well-known Hermitian and skew-Hermitian splitting (HSS) approach, with the particularity here of considering only one inner iteration. Experimental parallel performance comparison is conducted between the generalized minimal residual (GMRES) algorithm, the standard HSS and our asynchronous variant, on both real and complex non-Hermitian linear systems respectively arising from convection-diffusion and structural dynamics problems. A significant gain on execution time is observed in both cases.
\end{abstract}


\section{Introduction}

Many applications in scientific computing and engineering lead to the following system of linear equations,
\begin{equation}\label{eq:ls}
Ax=b,\quad A\in\mathbb{C}^{n\times n},\quad b\in\mathbb{C}^n.
\end{equation}
Let $A = M - N$ and $A = F - G$ be two splittings of $A$ with $M$ and $F$ being nonsingular.
The alternating iterative scheme for solving~\eqref{eq:ls} is defined as follows,
\begin{equation}\label{eq:alter}
\left\{\begin{array}{lcl}
Mx^{k+\frac{1}{2}} & = & Nx^{k} + b, \\
Fx^{k+1} & = & Gx^{k+\frac{1}{2}} + b,
\end{array}\right.
\end{equation}
which can be viewed as a stationary iterative scheme with an iteration matrix~$F^{-1}GM^{-1}N$.
Well-known early examples include the symmetric successive over-relaxation (SSOR) method~\cite{Sheldon1955,Conrad1979} and the alternating direction implicit (ADI) methods~\cite{Peaceman1955,Douglas1955,Marchuk1990}.
In~\cite{Benzi1997} the convergence of some alternating iterations were analyzed by eliminating the intermediate solution term~$x^{k+1/2}$ from~\eqref{eq:alter}; see also~\cite{Bai2003b}.
Recently, there has been growing interest in studies of the Hermitian and skew-Hermitian splitting (HSS) method~\cite{Bai2003} for solving~\eqref{eq:ls} when $A$ is non-Hermitian.
Let $\alpha>0$ be a given constant.
The HSS method can be written in the form
\begin{equation}\label{eq:hss}
\left\{\begin{array}{lcl}
(\alpha I + H)x^{k+\frac{1}{2}} & = & (\alpha I - S)x^{k} + b,\\
(\alpha I + S)x^{k+1} & = & (\alpha I - H)x^{k+\frac{1}{2}} + b,
\end{array}\right.
\end{equation}
where $H=(A+A^\mathsf{H})/2$ and $S=(A-A^\mathsf{H})/2$ are the Hermitian and skew-Hermitian parts of~$A$, respectively, and $I$ is the identity matrix.
Here, $A^\mathsf{H}$ denotes the conjugate transpose of~$A$.
This method can be obtained from~\eqref{eq:alter} by defining
\begin{equation}
\begin{array}{lcl}
M & := &\alpha I + H,\\
F & := &\alpha I + S.
\end{array}
\end{equation}
It was proved in~\cite{Bai2003} that when $H$ is positive definite, namely, $A$ is non-Hermitian positive definite, HSS converges unconditionally to the unique solution $x^{*}$ for any initial guess $x^{0}$.
The linear subsystems, however, especially the one involving~$\alpha I + S$, may still be difficult to solve, therefore much attention has been devoted to the inexact implementation.
More precisely, the tolerances for the inner iterative solvers may be relatively relaxed, while good convergence properties can still be retained according to numerical experiments; see~\cite{Bai2003,Benzi2007,Benzi2009,Bai2015}.
The HSS iterative scheme has been generalized to other splitting methods, as well as their preconditioned variants, for handling various problems in scientific computing; see,~e.g.,~\cite{Bertaccini2005,Li2007,Benzi2009,Bai2010,Wu2015,Li2015,Bai2019}.
There is also a number of studies on the optimal selection of~$\alpha$; see~\cite{Bai2003,Bai2006,Huang2014,Zou2020}.
The iterative scheme~\eqref{eq:hss} can be equivalently written in a residual-updating form, which achieves a higher accuracy at the cost of more computational effort; see~\cite{Bai2015} for a detailed discussion.

Parallel computing could be extremely useful when $A$ has large dimension.
In practice, the high cost of synchronization relative to that of computation is currently the major bottleneck in high-performance distributed computing systems, which motivates redesigning of parallel iterative algorithms.
One of the most interesting approaches, arising from basic relaxation methods, is the so-called asynchronous iterations~\cite{ChazMir1969,BertTsit1989}.
Asynchronous iterative scheme gives a full overlapping of communication and computation.
Every process has the flexibility to work at their own pace without waiting for the data acquisition.
A major difference between synchronous and asynchronous iterations lies in their predictability properties.
The former produces deterministic sequence of iterations, while the latter enables nondeterministic behaviors.
In~\cite{ChazMir1969} the first convergence result was established for the solution of linear systems, which was followed by the investigation of general fixed-point iterative models; see~\cite{Miellou1975,Baudet1978,ElTarazi1982,Bertsekas1983}.
In recent years, with the advent of very high-performance computing environment, asynchronous iterative scheme has gained much popularity.
The study of asynchronous domain decomposition methods, in both time and space domains, becomes an increasingly active area of research; see,~e.g.,~\cite{Magoules2017,Magoules2018,Magoules2018c,Magoules2018d,Yamazaki2019,ElHaddad2020}.
Another area that has seen growth in the last decades is the asynchronous convergence detection; see~\cite{Magoules2018e,GBenissan2020} and the references therein.

In this paper we focus on the asynchronous formulation of alternating iterations. In Section~\ref{sec:pre}, we recall some general tools and the asynchronous iterations theory used for the formulation and the convergence analysis of our asynchronous alternating scheme.
Section~\ref{sec:conv} presents the main contribution where we formulate our asynchronous alternating scheme and sufficient conditions for its convergence. Section~\ref{sec:exp} is devoted to numerical experiments on a parallel computing platform, featuring both a real three dimensional convection-diffusion problem and a complex two dimensional structural dynamic problem.
Finally, Section~\ref{sec:con} gives our conclusions.

\section{Generalities}\label{sec:pre}

\subsection{$\mathsf{H}$-matrix and $\mathsf{H}$-splitting}

In a general manner, let $\mathcal A_{i,j}$ denote the entry of a matrix $\mathcal A$ on its $i$-th row and $j$-th column, and let $x_{i}$ denote the $i$-th entry of a vector $x$. Comparisons $<$, $\le$, $>$, $\ge$ and $=$ between two matrices or vectors (of same shapes) are entrywise. The absolute value (or module) $|\mathcal A|$ of a matrix or a vector $\mathcal A$ is entrywise. The spectral radius of a matrix $\mathcal A$ is designated by $\rho(\mathcal A)$. In expressions like $\mathcal A < 0$ and like $x < 0$ with $\mathcal A$ and $x$ being a matrix and a vector, respectively, $0$ indicates a matrix and a vector, respectively, with all entries being $0$. $I$ stands for the identity matrix.

We recall now few general tools later used for the convergence analysis of the proposed asynchronous iterative method.
\begin{definition}
\label{def:mmatrix}
A square matrix $\mathcal A$ is an $\mathsf{M}$-matrix if and only if
\[
\exists \ \alpha \in \mathbb{R}: \quad \alpha I - \mathcal A \ge 0, \quad \alpha > \rho(\alpha I - \mathcal A).
\]
\end{definition}
\begin{definition}
\label{def:compmatrix}
The comparison matrix $\langle \mathcal A \rangle$ of a matrix $\mathcal A$  is defined as
\[
\langle \mathcal A \rangle_{i,i} := |\mathcal A_{i,i}|,
\qquad
\langle \mathcal A \rangle_{i,j} := -|\mathcal A_{i,j}|, \quad i \ne j.
\]
\end{definition}
\begin{definition}
\label{def:hmatrix}
A square matrix $\mathcal A$ is an $\mathsf{H}$-matrix if and only if its comparison matrix $\langle \mathcal A \rangle$ is an $\mathsf{M}$-matrix.
\end{definition}
\begin{lemma}
\label{lem:hmatrix}
A square matrix $\mathcal A$ is an $\mathsf{H}$-matrix if and only if
\[
\exists\ u > 0 : \quad \forall i, \ |\mathcal A_{i,i}| u_{i} > \sum_{j \ne i} |\mathcal A_{i,j}| u_{j}.
\]
\end{lemma}
\begin{proof}
This is directly implied by
Theorem 5' in \cite{Fan1958}.
\end{proof}

A splitting $\mathcal A = \mathcal M - \mathcal N$ of a matrix $\mathcal A$ consists of identifying a nonsingular matrix $\mathcal M$ and the resulting matrix $\mathcal N = \mathcal M - \mathcal A$, so as to define a relaxation operator
$
\mathcal M^{-1} \mathcal N = I - \mathcal M^{-1} \mathcal A.
$
\begin{definition}
\label{def:hsplit}
A splitting $\mathcal A = \mathcal M - \mathcal N$ is an $\mathsf H$-splitting if and only if $\langle \mathcal M \rangle - |\mathcal N|$ is an $\mathsf M$-matrix.
\end{definition}
\begin{lemma}
\label{lem:fs1992}
Let $\mathcal A = \mathcal M - \mathcal N$ be an $\mathsf H$-splitting. Then, we have
$
\rho(|I - \mathcal M^{-1} \mathcal A|) < 1.
$
\end{lemma}
\begin{proof}
This directly follows from Proof of Theorem 3.4 (c) in \cite{FromSzyld1992}.
\end{proof}
\begin{lemma}[refer to, e.g., Corollary 6.1 in \cite{BertTsit1989}]
\label{lem:pfc}
Let $\mathcal A$ be a square matrix. Then, we have
\[
\rho(\left|\mathcal A\right|) < 1 \quad \iff \quad \exists \ w > 0 : \ \left\|\mathcal A\right\|_{\infty}^{w} < 1,
\qquad\quad
\|\mathcal A\|_{\infty}^{w} := \max_{i} \frac{1}{w_{i}} \sum_{j} |\mathcal A_{i,j}| w_{j}.
\]
\end{lemma}

\subsection{Asynchronous iterations}

Consider, again, the linear system \eqref{eq:ls}, a splitting $A = M - N$ of the matrix $A$ and the resulting iterative scheme
\[
x^{k+1} = \left(I - M^{-1} A\right) x^{k} + M^{-1} b = x^{k} + M^{-1} \left(b - A x^{k}\right).
\]
Assume a distribution
\[
A =
\begin{bmatrix}
A^{(1)}\\
A^{(2)}\\
\vdots\\
A^{(m)}
\end{bmatrix},
\ \ 
b =
\begin{bmatrix}
b^{(1)}\\
b^{(2)}\\
\vdots\\
b^{(m)}
\end{bmatrix},
\ \ 
M =
\begin{bmatrix}
M^{(1)} & 0 & \cdots & 0\\
0 & M^{(2)} & \ddots & \vdots\\
\vdots & \ddots & \ddots & 0\\
0 & \cdots & 0 & M^{(m)}
\end{bmatrix}
\]
of both the system and the splitting of $A$. Note that the problem \eqref{eq:ls} can also corresponds to an augmented system resulting from a domain decomposition with overlapping subdomains, i.e., some rows in a submatrix $A^{(s_1)}$ are possibly replicated in another submatrix $A^{(s_2)}$, $s_1, s_2 \in \{1, \ldots, m\}$. A classical parallel relaxation is then given by
\[
\begin{aligned}
x^{(s),k+1} & = x^{(s),k} + {M^{(s)}}^{-1} \left(b^{(s)} - A^{(s)} \begin{bmatrix}x^{(1),k} & \cdots & x^{(m),k}\end{bmatrix}^{\mathsf T}\right) \quad \forall s \in \{1, \ldots, m\},\\
& = x^{(s),k} + {M^{(s)}}^{-1} \left(b^{(s)} - \sum_{q=1}^{m} A^{(s,q)} x^{(q),k}\right) \quad \forall s \in \{1, \ldots, m\}
\end{aligned}
\]
with
$
A^{(s)} = \begin{bmatrix}A^{(s,1)} & \cdots & A^{(s,m)}\end{bmatrix}.
$
The first feature of asynchronous iterations is the free steering (see, e.g., \cite{Schechter1959}), where, at each iteration $k$, a random subset $\Omega_{k} \subset \{1, \ldots, m\}$ of block-components can be updated. It is convenient to state a natural assumption,
\[
\operatorname{card} \left\{k \in \mathbb{N} : s \in \Omega_{k}\right\} = \infty \qquad \forall s \in \{1, \ldots, m\},
\]
which is implemented by the fact that no block-component stops being updated until convergence is globally reached. The second feature consists of modeling communication delays implying that at an iteration $k+1$, a block-component $s_1 \in \Omega_{k}$ is possibly updated using a block-component $s_2 \in \{1, \ldots, m\}$ computed at a random previous iteration $\delta_{s_1}(s_2, k) \le k$. It yields the parallel iterative scheme
\begin{equation}
\label{eq:ai}
x^{(s),k+1} = \left \{
\begin{array}{ll}
x^{(s),\delta_s(s,k)} + {M^{(s)}}^{-1} \left(b^{(s)} - \displaystyle\sum_{q=1}^{m} A^{(s,q)} x^{(q),\delta_s(q,k)}\right) & \forall s \in \Omega_{k},\\
x^{(s),k} & \forall s \notin \Omega_{k},
\end{array}
\right.
\end{equation}
where, as well, another natural assumption is made, stating that
\[
\lim_{k \to \infty} \delta_{s_1}(s_2,k) = \infty \qquad \forall s_1, s_2 \in \{1, \ldots, m\}.
\]
\begin{theorem}[Chazan and Miranker (1969) \cite{ChazMir1969}]
\label{theo:cm1969}
An asynchronous iterative method \eqref{eq:ai} converges from any initial guess $x^{0}$, with any sequence $\{\Omega_{k}\}_{k \in \mathbb{N}}$ and any functions $\delta_{1}$ to $\delta_{m}$ if and only if
$
\rho(|I - M^{-1} A|) < 1.
$
\end{theorem}
The model \eqref{eq:ai} was later generalized by Baudet \cite{Baudet1978} to arbitrary fixed-point iterations
\begin{equation}
\label{eq:gai}
x^{(s),k+1} = \left \{
\begin{array}{ll}
f^{(s)}\left(x^{(1),\delta_{s,1}(1,k)}, \ldots, x^{(m),\delta_{s,1}(m,k)},\right.&\\\left.\qquad\quad \ldots, x^{(1),\delta_{s,p}(1,k)}, \ldots, x^{(m),\delta_{s,p}(m,k)}\right) & \forall s \in \Omega_{k},\\
x^{(s),k} & \forall s \notin \Omega_{k},
\end{array}
\right.
\end{equation}
where the update of a block-component $s \in \Omega_{k}$ at an iteration $k$ depends on $p \in \mathbb{N}$ versions, $\delta_{s,1}(q,k)$ to $\delta_{s,p}(q,k)$, of each block-component $q \in \{1, \ldots, m\}$.
Let us denote by $\max (x, y)$ the vector given by
\[
(\max (x, y))_{i} := \max \{x_{i}, y_{i}\}
\]
with $x$ and $y$ being two vectors of same size. Let $X := (X_{1}, \ldots, X_{p})$ and $Y := (Y_{1}, \ldots, Y_{p})$ denote collections of $p$ vectors, i.e.,
\[
X_{t}^{} =
\begin{bmatrix}
X_{t}^{(1)} & \cdots & X_{t}^{(m)}
\end{bmatrix}^{\mathsf T},
\quad
Y_{t}^{} =
\begin{bmatrix}
Y_{t}^{(1)} & \cdots & Y_{t}^{(m)}
\end{bmatrix}^{\mathsf T},
\qquad
t \in \{1, \ldots, p\}.
\]
\begin{theorem}[Baudet (1978) \cite{Baudet1978}]
\label{theo:b1978}
An asynchronous iterative method \eqref{eq:gai} converges from any initial guess $x^{0}$, with any sequence $\{\Omega_{k}\}_{k \in \mathbb{N}}$ and any functions $\delta_{1,1}$ to $\delta_{m,p}$ if there exists a square matrix $\mathcal P$ such that $\mathcal P \ge 0$, $\rho(\mathcal P) < 1$ and
\[
\forall X, Y, \quad \left|f(X) - f(Y)\right| \le \mathcal P \max \left(\left|X_{1} - Y_{1}\right|, \ldots, \left|X_{p} - Y_{p}\right|\right).
\]
\end{theorem}

\section{Asynchronous alternating iterations}\label{sec:conv}

\subsection{Computational scheme}

Consider, now, the alternating scheme \eqref{eq:alter} which results in
\[
\begin{aligned}
x^{k+1} & = \left(I - F^{-1} A\right) x^{k+\frac{1}{2}} + F^{-1} b\\
& = \left(I - F^{-1} A\right) \left(I - M^{-1} A\right) x^{k} + \left(I - F^{-1} A\right) M^{-1} b + F^{-1} b\\
& = \left(I - F^{-1} \left(M + F - A\right) M^{-1} A\right) x^{k} + F^{-1} \left(M + F - A\right) M^{-1} b.
\end{aligned}
\]
Then, according to Theorem \ref{theo:cm1969}, such an induced parallel scheme is asynchronously convergent if
$
\rho\left(\left|I - F^{-1} \left(M + F - A\right) M^{-1} A\right|\right) < 1,
$
which is shown, in the next section, to be achieved under usual convergence conditions on the splittings $A = M - N$ and $A = F - G$.
Nevertheless, asynchronous relaxation based on such an operator cannot be implemented using the alternating form \eqref{eq:alter}, since the said operator is induced by strictly synchronizing $x^{k+\frac{1}{2}}$ and $x^{k+1}$.

Consider, then, an equivalent formulation of the alternating scheme \eqref{eq:alter},
\[
\left\{
\begin{array}{lcl}
y^{k} & := & x^{k} + M^{-1} \left(b - A x^{k}\right),\\
x^{k+1} & = & y^{k} + F^{-1} \left(b - A y^{k}\right),
\end{array}
\right.
\]
and assume that $F$ is distributed as $M$, i.e.,
\[
F =
\begin{bmatrix}
F^{(1)} & 0 & \cdots & 0\\
0 & F^{(2)} & \ddots & \vdots\\
\vdots & \ddots & \ddots & 0\\
0 & \cdots & 0 & F^{(m)}
\end{bmatrix}.
\]
Parallel asynchronous alternating methods are thus given by the computational scheme
\begin{equation}
\label{eq:async:alter}
\left \{
\begin{array}{lcl}
y^{(s),k} & := & x^{(s),\delta_s(s,k)}\\&&\quad +\ {M^{(s)}}^{-1} \left(b^{(s)} - \displaystyle\sum_{q=1}^{m} A^{(s,q)} x^{(q),\delta_s(q,k)}\right) \quad \forall s \in \{1, \ldots, m\},\\
x^{(s),k+1} & = & \left \{
\begin{array}{ll}
y^{(s),\delta_s(s,k)}&\\\quad +\ {F^{(s)}}^{-1} \left(b^{(s)} - \displaystyle\sum_{q=1}^{m} A^{(s,q)} y^{(q),\delta_s(q,k)}\right) & \forall s \in \Omega_{k},\\
x^{(s),k} & \forall s \notin \Omega_{k}.
\end{array}
\right.
\end{array}
\right.
\end{equation}
Assuming that the identity matrix $I$ is distributed as $A$, i.e.,
\[
I =
\begin{bmatrix}
I^{(1,1)} & \cdots & I^{(1,m)}\\
\vdots & \ddots & \vdots\\
I^{(m,1)} & \cdots & I^{(m,m)}
\end{bmatrix},
\]
it yields
\[
\begin{aligned}
x^{(s),k+1} & = \displaystyle\sum_{q=1}^{m} \left(I^{(s,q)} - {F^{(s)}}^{-1} A^{(s,q)}\right) y^{(q),\delta_s(q,k)} + {F^{(s)}}^{-1} b^{(s)}\\
& = \displaystyle\sum_{q=1}^{m} \left(I^{(s,q)} - {F^{(s)}}^{-1} A^{(s,q)}\right) \left(\displaystyle\sum_{r=1}^{m} \left(I^{(q,r)} - {M^{(q)}}^{-1} A^{(q,r)}\right) x^{(r),\delta_q(r,\delta_s(q,k))}\right.\\&\left.\qquad\qquad\qquad\qquad\qquad\qquad\qquad +\  {M^{(q)}}^{-1} b^{(q)}\right) + {F^{(s)}}^{-1} b^{(s)},\\
\end{aligned}
\]
which actually lies in the framework of the generalized model \eqref{eq:gai} with, here, $p = m$, since each update of a block-component depends on $m$ versions of the other block-components. Considering, then, a collection $X = \left(X_{1}, \ldots, X_{m}\right)$ of $m$ vectors, the corresponding mapping $f$ is given by
\[
\begin{aligned}
f^{(s)}(X) & := \displaystyle\sum_{q=1}^{m} \left(I^{(s,q)} - {F^{(s)}}^{-1} A^{(s,q)}\right) \left(\displaystyle\sum_{r=1}^{m} \left(I^{(q,r)} - {M^{(q)}}^{-1} A^{(q,r)}\right) X_{q}^{(r)}\right.\\&\left.\qquad\qquad\qquad\qquad\qquad\qquad\qquad
+\  {M^{(q)}}^{-1} b^{(q)}\right) + {F^{(s)}}^{-1} b^{(s)}\\
& = \displaystyle\sum_{q=1}^{m} P_{q}^{(s)} X_{q} + \left(I^{(s)} - {F^{(s)}}^{-1} A^{(s)}\right) M^{-1} b + {F^{(s)}}^{-1} b^{(s)},\\
f(X) & := \displaystyle\sum_{q=1}^{m} P_{q} X_{q} + \left(I - F^{-1} A\right) M^{-1} b + F^{-1} b
\end{aligned}
\]
with
$
P_{q}^{(s)} := \left(I^{(s,q)} - {F^{(s)}}^{-1} A^{(s,q)}\right) \left(I^{(q)} - {M^{(q)}}^{-1} A^{(q)}\right), \ q,s \in \{1, \ldots, m\},
$
and
$
P_{q} := \begin{bmatrix}P_{q}^{(1)} & \cdots & P_{q}^{(m)}\end{bmatrix}^{\mathsf T}, \ q \in \{1, \ldots, m\}.
$

\subsection{Convergence conditions}

We analyze, now, sufficient conditions for the convergence of our asynchronous alternating iterative scheme \eqref{eq:async:alter}. To the best of our knowledge, Lemma~\ref{lem:cwmwmn}, Proposition~\ref{prop:ahss} and Corollary~\ref{cor:ahss} are new. Proposition~\ref{prop:ahss} and Corollary~\ref{cor:ahss} highlight how combining properties of the operators $I - F^{-1} A$ and $I - M^{-1} A$ imply a resulting contracting operator $\left(I - F^{-1} A\right)\left(I - M^{-1} A\right)$. Our main results consist of Theorem~\ref{theo:ahss_practical} and Corollary~\ref{cor:ahss_practical} where the same combined conditions are shown to be sufficient for the convergence of asynchronous alternating methods \eqref{eq:async:alter}, despite the induced, slightly different, iterations operator.

Let, first, $\mathcal A$ be a matrix with arbitrary shape, let $w$ be a vector with as many entries as the number of columns in $\mathcal A$, and let $v$ be a vector with as many entries as the number of rows in $\mathcal A$, and with no $0$ entry. Let $\tau(\mathcal A, w, v)$ denote the vector given by the row-sums
\[
\tau_{i}(\mathcal A, w, v) := \left(\tau(\mathcal A, w, v)\right)_{i} := \frac{1}{v_{i}} \sum_{j} \left|\mathcal A_{i,j}\right| w_{j} \qquad \forall i.
\]
Note, then, that, for a square matrix $\mathcal A$,
\[
\|\mathcal A\|_{\infty}^{w} = \max_{i} \tau_{i}(\mathcal A, w, w), \qquad w > 0.
\]
\begin{lemma}
\label{lem:cwmwmn}
Let $\mathcal A$ and $\mathcal B$ be matrices with shapes such that $\mathcal A \mathcal B$ is calculable. Let $u > 0$, $v > 0$ and $w$ be vectors with dimensions such that $\tau(\mathcal A, u, v)$ and $\tau(\mathcal B, w, u)$ are calculable. Then, we have
\[
\tau(\mathcal B, w, u) <
\begin{bmatrix}
1 & 1 & \cdots & 1
\end{bmatrix}^{\mathsf T}
\quad \implies \quad \tau(\mathcal A \mathcal B, w, v) < \tau(\mathcal A, u, v).
\]
\end{lemma}
\begin{proof}
Let us index rows and columns of $\mathcal A$ by $i$ and $j$, respectively, and columns of $\mathcal B$ by $l$. We have
\[
\begin{aligned}
\tau_{i}(\mathcal A \mathcal B, w, v) := \frac{1}{v_{i}} \sum_{l} \left|(\mathcal A \mathcal B)_{i,l}\right| w_{l} & = \frac{1}{v_{i}} \sum_{l} \left|\sum_{j} \mathcal A_{i,j} \mathcal B_{j,l}\right| w_{l}\\
& \le \frac{1}{v_{i}} \sum_{l} \sum_{j} \left|\mathcal A_{i,j} \mathcal B_{j,l}\right| w_{l}\\
& = \frac{1}{v_{i}} \sum_{l} \sum_{j} \frac{1}{u_{j}} \left|\mathcal A_{i,j}\right| \left|\mathcal B_{j,l}\right| u_{j} w_{l}\\
& = \frac{1}{v_{i}} \sum_{j} \left(\frac{1}{u_{j}} \sum_{l} \left|\mathcal B_{j,l}\right| w_{l}\right) \left|\mathcal A_{i,j}\right| u_{j}\\
& = \frac{1}{v_{i}} \sum_{j} \tau_{j}(\mathcal B, w, u) \left|\mathcal A_{i,j}\right| u_{j}.
\end{aligned}
\]
It yields that if $\tau_{j}(\mathcal B, w, u) < 1$ for all $j$, then
\[
\begin{aligned}
&& \tau_{j}(\mathcal B, w, u) \left|\mathcal A_{i,j}\right| u_{j} & < \left|\mathcal A_{i,j}\right| u_{j} \quad \forall j \ \forall i,\\
&& \frac{1}{v_{i}} \sum_{j} \tau_{j}(\mathcal B, w, u) \left|\mathcal A_{i,j}\right| u_{j} & < \frac{1}{v_{i}} \sum_{j} \left|\mathcal A_{i,j}\right| u_{j} \quad \forall i,\\
\frac{1}{v_{i}} \sum_{l} \left|(\mathcal A \mathcal B)_{i,l}\right| w_{l} & \ \le & \frac{1}{v_{i}} \sum_{j} \tau_{j}(\mathcal B, w, u) \left|\mathcal A_{i,j}\right| u_{j} & < \frac{1}{v_{i}} \sum_{j} \left|\mathcal A_{i,j}\right| u_{j} \quad \forall i,\\
\tau_{i}(\mathcal A \mathcal B, w, v) &&& < \tau_{i}(\mathcal A, u, v) \quad \forall i,
\end{aligned}
\]
which concludes the proof.
\end{proof}
\begin{proposition}
\label{prop:ahss}
Let
\[
Q :=
\begin{bmatrix}
0 & I - M^{-1} A\\
I - F^{-1} A & 0
\end{bmatrix}.
\]
We have
\[
\rho(|Q|) < 1 \quad \implies \quad \rho\left(\left|I - F^{-1} \left(M + F - A\right) M^{-1} A\right|\right) < 1.
\]
\end{proposition}
\begin{proof}
According to Lemma \ref{lem:pfc},
\[
\rho(|Q|) < 1 \quad \iff \quad \exists \ W > 0 : \ \|Q\|_{\infty}^{W} < 1.
\]
According to the two blocks of $Q$, take
$
W =
\begin{bmatrix}
W_{1} & W_{2}
\end{bmatrix}^{\mathsf T}.
$
Then, we have both
\[
\left\{
\begin{array}{lcl}
\tau\left(I - M^{-1} A, W_{2}, W_{1}\right) & < & \begin{bmatrix} 1 & 1 & \cdots & 1 \end{bmatrix}^{\mathsf T},\\
\tau\left(I - F^{-1} A, W_{1}, W_{2}\right) & < & \begin{bmatrix} 1 & 1 & \cdots & 1 \end{bmatrix}^{\mathsf T}.
\end{array}
\right.
\]
Lemma \ref{lem:cwmwmn} therefore ensures
\[
\begin{aligned}
\tau\left(\left(I - F^{-1} A\right) \left(I - M^{-1} A\right), W_{2}, W_{2}\right) & < \tau\left(I - F^{-1} A, W_{1}, W_{2}\right)\\
& < \begin{bmatrix} 1 & 1 & \cdots & 1 \end{bmatrix}^{\mathsf T},
\end{aligned}
\]
which leads to
$
\left\|\left(I - F^{-1} A\right) \left(I - M^{-1} A\right)\right\|_{\infty}^{W_{2}} < 1.
$
Recall that
\[
\left(I - F^{-1} A\right) \left(I - M^{-1} A\right) = I - F^{-1} \left(M + F - A\right) M^{-1} A.
\]
Lemma \ref{lem:pfc} finally ensures
$
\rho\left(\left|I - F^{-1} \left(M + F - A\right) M^{-1} A\right|\right) < 1,
$
which concludes the proof.
\end{proof}
\begin{corollary}
\label{cor:ahss}
if $A$ is an $\mathsf H$-matrix, then
\[
\left\{
\begin{array}{lcl}
\langle M \rangle - |M - A| & = & \langle A \rangle,\\
\langle F \rangle - |F - A| & = & \langle A \rangle
\end{array}
\right.
\quad \implies \quad \rho\left(\left|I - F^{-1} \left(M + F - A\right) M^{-1} A\right|\right) < 1.
\]
\end{corollary}
\begin{proof}
Considering that $A$ is an $\mathsf H$-matrix, take $u > 0$ like in Lemma \ref{lem:hmatrix}, so as to have
\[
|A_{i,i}| u_{i} > \sum_{j \ne i} |A_{i,j}| u_{j} \quad \forall i.
\]
We also have
\[
\langle M \rangle - |M - A| = \langle A \rangle \quad \implies \quad \forall i, \ 
\left\{
\begin{array}{lcl}
|M_{i,i}| - |M_{i,i} - A_{i,i}| & = & |A_{i,i}|,\\
- |M_{i,j}| - |M_{i,j} - A_{i,j}| & = & - |A_{i,j}| \quad \forall j \ne i,
\end{array}
\right.
\]
and, then,
\[
\left\{
\begin{array}{lcl}
|M_{i,i}| u_{i} - |M_{i,i} - A_{i,i}| u_{i} & = & |A_{i,i}| u_{i},\\
- |M_{i,j}| u_{j} - |M_{i,j} - A_{i,j}| u_{j} & = & - |A_{i,j}| u_{j} \quad \forall j \ne i.
\end{array}
\right.
\]
It yields that, $\forall i$,
\[
\begin{aligned}
|M_{i,i}| u_{i} - \sum_{j \ne i} |M_{i,j}| u_{j} - |M_{i,i} - A_{i,i}| u_{i} - \sum_{j \ne i} |M_{i,j} - A_{i,j}| u_{j} & = |A_{i,i}| u_{i} - \sum_{j \ne i} |A_{i,j}| u_{j}\\
& > 0,
\end{aligned}
\]
which implies, with $F$ also satisfying $\langle F \rangle - |F - A| = \langle A \rangle$, that the matrix
\[
\widehat A :=
\begin{bmatrix}
M & A - M\\
A - F & F
\end{bmatrix}
\]
is an $\mathsf H$-matrix, according to Lemma \ref{lem:hmatrix}. Define, then,
\[
\widehat M :=
\begin{bmatrix}
M & 0\\
0 & F
\end{bmatrix},
\]
and note that $\left\langle \widehat M \right\rangle - \left|\widehat M - \widehat A\right| = \left\langle \widehat A \right\rangle$, which implies, by Definition \ref{def:hmatrix}, that $\left\langle \widehat M \right\rangle - \left|\widehat M - \widehat A\right|$ is an $\mathsf M$-matrix, hence, by Definition \ref{def:hsplit}, $\widehat A = \widehat M - \left(\widehat M - \widehat A\right)$ is an $\mathsf H$-splitting. Lemma \ref{lem:fs1992} therefore ensures that
$
\rho\left(\left|\widehat M^{-1} \left(\widehat M - \widehat A\right)\right|\right) < 1,
$
and one can verify that
\[
\widehat M^{-1} \left(\widehat M - \widehat A\right) =
\begin{bmatrix}
0 & I-M^{-1}A\\
I-F^{-1}A & 0
\end{bmatrix}.
\]
Proposition \ref{prop:ahss} therefore finally applies, which concludes the proof.
\end{proof}
\begin{theorem}
\label{theo:ahss_practical}
Let
\[
Q :=
\begin{bmatrix}
0 & I - M^{-1} A\\
I - F^{-1} A & 0
\end{bmatrix}.
\]
An asynchronous alternating method \eqref{eq:async:alter} converges from any initial guess $x^{0}$, with any sequence $\{\Omega_{k}\}_{k \in \mathbb{N}}$ and any functions $\delta_{1}$ to $\delta_{m}$ if $\rho(|Q|) < 1$.
\end{theorem}
\begin{proof}
Consider two collections, $X = \left(X_1, \ldots, X_m\right)$ and $Y = \left(Y_1, \ldots, Y_m\right)$, of $m$ vectors. We have
\[
\begin{aligned}
|f(X) - f(Y)| & = \left|\sum_{q=1}^{m} P_{q} \left(X_{q} - Y_{q}\right)\right|\\
& \le \sum_{q=1}^{m} \left|P_{q}\right| \max \left(\left|X_{1} - Y_{1}\right|, \ldots, \left|X_{m} - Y_{m}\right|\right).
\end{aligned}
\]
Consequently, according to Theorem \ref{theo:b1978}, an asynchronous alternating method \eqref{eq:async:alter} is convergent if
$
\rho\left(\sum_{q=1}^{m} \left|P_{q}\right|\right) < 1.
$
Recall, then, that according to Lemma \ref{lem:pfc},
\[
\rho(|Q|) < 1 \quad \iff \quad \exists \ W > 0 : \ \|Q\|_{\infty}^{W} < 1.
\]
According to the two blocks of $Q$, take
$
W =
\begin{bmatrix}
W_{1} & W_{2}
\end{bmatrix}^{\mathsf T}.
$
Then, we have both
\[
\left\{
\begin{array}{lcl}
\tau\left(I - M^{-1} A, W_{2}, W_{1}\right) & < & \begin{bmatrix} 1 & 1 & \cdots & 1 \end{bmatrix}^{\mathsf T},\\
\tau\left(I - F^{-1} A, W_{1}, W_{2}\right) & < & \begin{bmatrix} 1 & 1 & \cdots & 1 \end{bmatrix}^{\mathsf T},
\end{array}
\right.
\]
implying, as well,
\[
\tau\left(I^{(q)} - {M^{(q)}}^{-1} A^{(q)}, W_{2}, W_{1}^{(q)}\right) < \begin{bmatrix} 1 & 1 & \cdots & 1 \end{bmatrix}^{\mathsf T} \quad \forall q \in \{1, \ldots, m\}.
\]
Lemma \ref{lem:cwmwmn} therefore ensures, with $s \in \{1, \ldots, m\}$,
\[
\begin{aligned}
\tau\left(\left(I^{(s,q)} - {F^{(s)}}^{-1} A^{(s,q)}\right) \left(I^{(q)} - {M^{(q)}}^{-1} A^{(q)}\right), W_{2}, W_{2}^{(s)}\right) & < \tau\left(I^{(s,q)} - {F^{(s)}}^{-1} A^{(s,q)},\right.\\&\left.\qquad\quad W_{1}^{(q)}, W_{2}^{(s)}\right).
\end{aligned}
\]
Recall that
$
P_{q}^{(s)} := \left(I^{(s,q)} - {F^{(s)}}^{-1} A^{(s,q)}\right) \left(I^{(q)} - {M^{(q)}}^{-1} A^{(q)}\right), \ q,s \in \{1, \ldots, m\}.
$
Then, we have
\[
\begin{aligned}
\tau\left(P^{(s)}_{q}, W_{2}, W_{2}^{(s)}\right) & < \tau\left(I^{(s,q)} - {F^{(s)}}^{-1} A^{(s,q)}, W_{1}^{(q)}, W_{2}^{(s)}\right),\\
\tau\left(\left|P^{(s)}_{q}\right|, W_{2}, W_{2}^{(s)}\right) & < \tau\left(I^{(s,q)} - {F^{(s)}}^{-1} A^{(s,q)}, W_{1}^{(q)}, W_{2}^{(s)}\right),\\
\sum_{q=1}^{m}\tau\left(\left|P^{(s)}_{q}\right|, W_{2}, W_{2}^{(s)}\right) & < \sum_{q=1}^{m}\tau\left(I^{(s,q)} - {F^{(s)}}^{-1} A^{(s,q)}, W_{1}^{(q)}, W_{2}^{(s)}\right),\\
\tau\left(\sum_{q=1}^{m}\left|P^{(s)}_{q}\right|, W_{2}, W_{2}^{(s)}\right) & < \tau\left(I^{(s)} - {F^{(s)}}^{-1} A^{(s)}, W_{1}^{}, W_{2}^{(s)}\right),\\
\tau\left(\sum_{q=1}^{m}\left|P^{}_{q}\right|, W_{2}, W_{2}^{}\right) & < \tau\left(I - F^{-1} A, W_{1}^{}, W_{2}^{}\right),\\
& < \begin{bmatrix} 1 & 1 & \cdots & 1 \end{bmatrix}^{\mathsf T},
\end{aligned}
\]
which leads to
$
\left\|\sum_{q=1}^{m}\left|P^{}_{q}\right|\right\|_{\infty}^{W_{2}} < 1.
$
By Lemma \ref{lem:pfc}, we therefore satisfy
$
\rho\left(\sum_{q=1}^{m}\left|P^{}_{q}\right|\right) < 1,
$
which concludes the proof.
\end{proof}
\begin{corollary}
\label{cor:ahss_practical}
An asynchronous alternating method \eqref{eq:async:alter} converges from any initial guess $x^{0}$, with any sequence $\{\Omega_{k}\}_{k \in \mathbb{N}}$ and any functions $\delta_{1}$ to $\delta_{m}$ if $A$ is an $\mathsf H$-matrix and
\[
\left\{
\begin{array}{lcl}
\langle M \rangle - |M - A| & = & \langle A \rangle,\\
\langle F \rangle - |F - A| & = & \langle A \rangle.
\end{array}
\right.
\]
\end{corollary}
\begin{proof}
This follows in the same way as Corollary \ref{cor:ahss}.
\end{proof}
Let $\mathcal D(\mathcal A)$ denote the diagonal matrix obtained from the diagonal of a matrix $\mathcal A$.
\begin{remark}
For practical applications of Corollary \ref{cor:ahss_practical}, let $\Lambda$ be a diagonal real matrix such that
$
\Lambda_{i,i} \ge 1 \ \forall i.
$
We straightforwardly have
\[
\mathcal M = \Lambda \mathcal D(\mathcal A) \quad \implies \quad \langle \mathcal M \rangle - |\mathcal M - \mathcal A| = \langle \mathcal A \rangle.
\]
\end{remark}
\begin{remark}
In regard to the HSS splitting, if $A$ is a real matrix with $\mathcal D(A) \ge 0$, and splitting matrices $M$ and $F$ are given by
\[
M := \mathcal D(\alpha I + H), \qquad F := \mathcal D(\alpha I + S), \qquad \alpha \ge \max_{i} A_{i,i},
\]
then we have both
\[
M = \alpha I + \mathcal D(A) \ge \mathcal D(A), \qquad F = \alpha I \ge \mathcal D(A),
\]
which satisfy
$
M = \Lambda_{M} \mathcal D(A), \ F = \Lambda_{F} \mathcal D(A),
$
where $\Lambda_{M}$ and $\Lambda_{F}$ are two diagonal real matrices with entries greater than or equal to $1$.
\end{remark}

{

\section{Implementation aspects}\label{sec:impl}

The two alternating iterations of the HSS method require the solution of two secondary problems involving the coefficient matrices $\alpha I + H$ and $\alpha I + S$, respectively. In practice, as pointed out in, e.g., \cite{Bai2003, Wu2015}, these problems are inexactly solved by means of iterative algorithms. A general description for both HSS and inexact HSS (IHSS) can be given by Algorithm~\ref{algo:hss}.
\begin{algorithm}[htbp]
\caption{HSS(solverH, solverS)}
\label{algo:hss}
{
\begin{algorithmic}[1]
\STATE{$x$ := $x^{0}$}
\STATE{$r$ := $b - A x$}
\STATE{$k$ := $0$}
\WHILE{$\|r\| > \varepsilon \|b\|$ and $k < k_{\text{max}}$}
  \STATE{$y$ := solverH.solve($\alpha I + H$, $r$)}
  \STATE{$x$ := $x + y$}
  \STATE{$r$ := $b - A x$}
  \STATE{$y$ := solverS.solve($\alpha I + S$, $r$)}
  \STATE{$x$ := $x + y$}
  \STATE{$r$ := $b - A x$}
  \STATE{$k$ := $k + 1$}
\ENDWHILE
\end{algorithmic}
}
\end{algorithm}
We can then designate by, e.g, HSS(CG, GMRES) an IHSS algorithm with the conjugate gradient (CG) method~\cite{HestenesStiefel1952} for solving the shifted Hermitian problem and the generalized minimal residual (GMRES) method~\cite{SaadSchultz1986} for solving the shifted skew-Hermitian one.

Asynchronous HSS iterations necessarily belong to the class of IHSS algorithms since they obviously require the inner solvers to be asynchronous too, which further reduces such an approach to the subclass of IHSS with inner splittings. Taking, then, e.g., a splitting
$
\alpha I + H = M - N,
$
the solution, at each outer iteration $k$, of
\[
(\alpha I + H) y^{k} = b - A x^{k}
\]
can be given by several inner iterations
\begin{equation}
\label{eq:inner_iter}
y^{k,l+1} = y^{k,l} + M^{-1} (b - A x^{k} - (\alpha I + H) y^{k,l}),
\end{equation}
where $l$ is the inner iteration variable. Furthermore, when dealing with two-stage asynchronous iterations, one should particularly take advantage of the possibility to use the inner solution vector $y^{k,l+1}$ with any value of $l$, given that asynchronous relaxation is very likely to benefit from each newly updated data. We refer the reader to, e.g., \cite{ElBazEtAl1996, FromSzyld1998} for more insights into the so called ``asynchronous iterations with flexible communication''. Moreover, analysis of matrix splittings for two-stage asynchronous iterations reveals that convergence of such methods can be guaranteed for any number of inner iterations (see, e.g., \cite{FromSzyld1994}). According, therefore, to efficiency aspects related to flexible communication ideas, it is of some interest, in the end, to simply consider only one iteration of \eqref{eq:inner_iter}. If, in particular, we also consider as initial guess $y^{k,0} := 0$, then we can define
\[
y^{k} := y^{k,1} = M^{-1} (b - A x^{k}),
\]
so as to finally have
\[
x^{k+\frac{1}{2}} = x^{k} + M^{-1} (b - A x^{k}),
\]
which falls under the general alternating scheme \eqref{eq:alter} that has been considered in our theoretical analysis. Such a specialization of Algorithm~\ref{algo:hss} is given by Algorithm~\ref{algo:hss_inner1}, where $M^{-1}$ and $F^{-1}$ are preconditioners of $\alpha I + H$ and $\alpha I + S$, respectively.
\begin{algorithm}[htbp]
\caption{HSS($M^{-1}$, $F^{-1}$)}
\label{algo:hss_inner1}
{
\begin{algorithmic}[1]
\STATE{$x$ := $x^{0}$}
\STATE{$r$ := $b - A x$}
\STATE{$k$ := $0$}
\WHILE{$\|r\| > \varepsilon \|b\|$ and $k < k_{\text{max}}$}
  \STATE{$x$ := $x + M^{-1} r$}
  \STATE{$r$ := $b - A x$}
  \STATE{$x$ := $x + F^{-1} r$}
  \STATE{$r$ := $b - A x$}
  \STATE{$k$ := $k + 1$}
\ENDWHILE
\end{algorithmic}
}
\end{algorithm}
Note that Algorithm~\ref{algo:hss_inner1} needs to be specifically implemented instead of just using Algorithm~\ref{algo:hss} with calls of relaxation-based inner solvers with maximum number of iterations set to $1$. Indeed, on pure computer science aspects, avoiding inner function calls and loops can result in a very significant execution time saving, which even makes HSS($M^{-1}$, $F^{-1}$) possibly competitive, in practice, with, e.g., HSS(CG, GMRES), as we shall see in Section~\ref{sec:exp}.

From Algorithm~\ref{algo:hss_inner1}, iterative scheme \eqref{eq:async:alter}, programming models \cite{Magoules2017b,Magoules2018b} and convergence detection approach \cite{GBenissan2020}, asynchronous parallel implementation of HSS iterations is obtained as described by Algorithm~\ref{algo:hss_async}, where the communication routines start with ``Com'' and are blocking by default. Their non-blocking counterparts are designated by ``ICom'' with the letter ``I'' standing for ``immediate'', similarly to the Message Passing Interface (MPI) standard.
\begin{algorithm}[htbp]
\caption{Asynchronous parallel HSS(${M^{(s)}}^{-1}$, ${F^{(s)}}^{-1}$) on process $s \in \{1, \ldots, m\}$}
\label{algo:hss_async}
{
\begin{algorithmic}[1]
\STATE{$x^{(s)}$ := $x^{(s),0}$}
\STATE{$x$ := IComSendRecvInit($x^{(s)}$)}
\STATE{$r^{(s)}$ := $b^{(s)} - A^{(s)} x$}
\STATE{${rr}^{(s)}$ := ${r^{(s)}}^{\mathsf H} r^{(s)}$}
\STATE{$rr$ := ComSum(${rr}^{(s)}$)}
\STATE{$\|r\|$ := $\sqrt{rr}$}
\STATE{$\tau$ := False}
\STATE{$k$ := $0$}
\WHILE{$\|r\| > \varepsilon \|b\|$ and $k < k_{\text{max}}$}
  \STATE{$x^{(s)}$ := $x^{(s)} + {M^{(s)}}^{-1} r^{(s)}$}
  \STATE{$x$ := IComSendRecv($x^{(s)}$)}
  \STATE{$r^{(s)}$ := $b^{(s)} - A^{(s)} x$}
  \STATE{$x^{(s)}$ := $x^{(s)} + {F^{(s)}}^{-1} r^{(s)}$}
  \STATE{$x$ := IComSendRecv($x^{(s)}$)}
  \STATE{$r^{(s)}$ := $b^{(s)} - A^{(s)} x$}
  \IF{\NOT $\tau$}
    \STATE{${rr}^{(s)}$ := ${r^{(s)}}^{\mathsf H} r^{(s)}$}
    \STATE{ComRequest := IComSum(${rr}^{(s)}$, $rr$)}
    \STATE{$\tau$ := True}
  \ENDIF
  \STATE{$\sigma$ := ComTest(ComRequest)}
  \IF{$\sigma$}
    \STATE{$\|r\|$ := $\sqrt{rr}$}
    \STATE{$\tau$ := False}
    \STATE{$k$ := $k + 1$}
  \ENDIF
\ENDWHILE
\end{algorithmic}
}
\end{algorithm}
The routines ComSum and IComSum are used to compute dot product $r^{\mathsf H} r$ with $r = b - A x$ by global reduction operation
\[
\sum_{q=1}^{m} {r^{(q)}}^{\mathsf H} r^{(q)}, \qquad r^{(q)} = b^{(q)} - A^{(q)} x.
\]
They can readily be replaced by MPI routines MPI\_Allreduce and MPI\_Iallreduce, respectively. The object ComRequest and the routine ComTest are therefore analogous to MPI\_Request and MPI\_Test. Such a simple way to reliably use the classical loop stopping criterion $\|r\| > \varepsilon \|b\|$ in case of asynchronous iterations is due to \cite{GBenissan2020}. It also allows for considering a counter, $k$, of the number of global convergence tests. On the other hand, the data exchange routine IComSendRecv has to be a bit constructed using, e.g., MPI routines MPI\_Isend and MPI\_Irecv. Briefly, the routine IComSendRecvInit triggers non-blocking requests for message sending ($x^{(s)}$) and reception ($x^{(q)}$, $q \ne s$), and fills up the components $x^{(q)}$, $q \ne s$, of the vector $x$ with any arbitrary values. Note that both storage and communication of components $x^{(q)}$, $q \ne s$, should actually be limited to values which are necessary for computing the product $A^{(s)} x$, according to the nonzero entries in $A^{(s)}$. The subsequent calls to the routine IComSendRecv then check completion of previous requests, update $x$ with received data and trigger new instances of the completed requests. Further details can be found in, e.g., \cite{Magoules2018b}.

}

\section{Numerical experiments}\label{sec:exp}

{
\subsection{Problems and overall settings}
}

Numerical experiments have been conducted on two kinds of problem. The first one consists of a three-dimensional (3D) convection-diffusion equation,
\begin{equation}\label{eq:cd}
-\Delta u + c\cdot\nabla u = f\mbox{ in $\Omega$}
\end{equation}
with $\Omega=[0,1]\times[0,1]\times[0,1]$ and Dirichlet boundary conditions. Discretization has been achieved using seven-point centered differences for both convection and diffusion terms. A fixed value, $20$, has been used for all elements in the three-dimensional vector $c$ as convection parameter. The entries of the exact discrete solution, $x^*$, have been taken randomly in $[0,1)$ and the right-hand side has then been constructed as $b=Ax^*$.

The second kind of problem consists of a 2D structural dynamics equation (see, e.g., \cite{Benzi2008, Bai2010}),
\begin{equation}\label{eq:sd}
\left[\left(-\omega^2 L + K\right) + \operatorname{i} \left(\omega C_v + C_h\right)\right] x = b,
\end{equation}
where $L$ and $K$ denote the mass and stiffness matrices, respectively; $C_v$ and $C_h$ denote the viscous and hysteretic damping matrices, respectively; $\omega$ denotes the circular frequency. The values of the matrices and the parameters have been taken from \cite{Bai2010}. The matrix $K$ is the five-point finite difference discretization of a diffusion term on the unit square $[0,1]\times[0,1]$ with Dirichlet boundary conditions. The other matrices have been set as $L=I$, $C_v=10I$, $C_h=\mu K$, where $\mu=0.02$, and $I$ denotes the $n \times n$ identity matrix. The circular frequency $\omega$ has been set to $\pi$. The right-hand side has been taken as $b=(1+\operatorname{i})Aq$ with $q$ being a vector of $1$, to ensure that all entries of $x^*$ equal $1+\operatorname{i}$.

In the following, parallel execution times (wall-clock), numbers of iterations, $k$, and final residual errors, $r$, are reported for the GMRES~\cite{SaadSchultz1986}, the IHSS~\cite{Bai2003} (Algorithms \ref{algo:hss} and \ref{algo:hss_inner1}) and the asynchronous IHSS methods (Algorithm \ref{algo:hss_async}), with a stopping criterion set so as to have
\[
r = \frac{\|b-Ax^{*}\|}{\|b\|} < 10^{-6}.
\]
In case of asynchronous execution, minimum and maximum numbers of local iterations, $k_{\text{min}}$ and $k_{\text{max}}$, respectively, are considered since there is not global iterations $k$. Both for synchronous and asynchronous HSS(${M}^{-1}$, ${F}^{-1}$) (respectively, Algorithms \ref{algo:hss_inner1} and \ref{algo:hss_async}), we took
\[
M := \mathcal D(\alpha I + H), \qquad F := \mathcal D(\alpha I + S).
\]
All of the tests have been entirely implemented in the Python language, using NumPy, SciPy Sparse and MPI4Py~\cite{Dalcin2005} modules.

{

A comparison with some results in \cite{Bai2010} about the problem \eqref{eq:sd} (Example 4.2 in \cite{Bai2010}) is reported in Table~\ref{tab:comp_bai2010} for single-process execution of full GMRES, GMRES(restart), and HSS(CG, GMRES(restart)) with inner residual threshold set to $10^{-10}$ in order to compare with an ``exact'' HSS.
\begin{table}[htbp]
\caption{Comparison with Ref. \cite{Bai2010} for the test case \eqref{eq:sd}, number of processes $p = 1$.}
{
\begin{tabular}{cc}
\hline
Experiment & Results\\
\hline
\begin{tabular}{c}
Ref. \cite{Bai2010}\\MATLAB\\2.66 GHz CPU\\1.97 GB RAM
\end{tabular} &
\begin{tabular}{ccccc}
  $n$ & \multicolumn{2}{c}{$64 \times 64$} & \multicolumn{2}{c}{$128 \times 128$}\\
  \hline
  Method & Clock (sec) & $k$ & Clock (sec) & $k$\\
  \hline
  HSS & 4.81 & 284 & 60 & 540\\
  GMRES(10) & 1.08 & 973 & 20 & 3096\\
  GMRES(20) & 1.50 & 632 & 22 & 1704\\
  GMRES & 2.98 & 161 & 45 & 308
\end{tabular}\\
\hline
\begin{tabular}{c}
Python\\2.40 GHz CPU\\174 GB RAM
\end{tabular} &
\begin{tabular}{ccccc}
  $n$ & \multicolumn{2}{c}{$64 \times 64$} & \multicolumn{2}{c}{$128 \times 128$}\\
  \hline
  Method & Clock (sec) & $k$ & Clock (sec) & $k$\\
  \hline
  HSS(CG,GMRES(10)) & 4.80 & 284 & 44 & 540\\
  GMRES(10) & 0.36 & 1072 & 3.56 & 3346\\
  GMRES(20) & 0.33 & 672 & 2.70 & 1790\\
  GMRES & 0.44 & 161 & 5.19 & 308
\end{tabular}\\
\hline
\end{tabular}
}
\label{tab:comp_bai2010}
\end{table}
The experimentally optimal value of $\alpha$, according to \cite{Bai2010}, was considered for each problem size $n$ ($\alpha = 0.12$ for $n = 64^{2}$, and $\alpha = 0.07$ for $n = 128^{2}$). We recall that the experiments in \cite{Bai2010} were run in MATLAB on a personal computer consisting of a 2.66 GHz Intel Core Duo central processing unit (CPU) and 1.97 GB of random access memory (RAM). Our single-process tests, here, have been performed on a computational cluster node consisting of a 2.40 GHz Intel Xeon Skylake CPU and 174 GB of RAM. Same numbers of iterations are obtained for our implementation of HSS(CG, GMRES(10)), where both CG and GMRES's tolerances were set to $10^{-10}$, and the HSS experimented in \cite{Bai2010} with direct inner solvers. Same result is observed for full GMRES too, while very slight differences appear for the restarted GMRES.

The remaining tests, which involve multi-process execution, have been performed on cluster nodes consisting of 2 $\times$ 12-cores 2.30 GHz Intel Xeon Haswell CPU (24 cores per node) and
48 GB of RAM (2 GB per core).
The nodes are interconnected through a 56 Gb/s fourteen data rate (FDR) Infiniband network, on which the SGI MPT library is used as implementation of the MPI standard.

\subsection{Results on the 3D convection-diffusion problem}

\subsubsection{Optimal parameters}

The 3D convection-diffusion test case \eqref{eq:cd} was run on an obtained discrete problem with $n = 100^{3}$ unknowns, using from $p = 48$ to $p = 192$ processor cores (one MPI process per core).

Table \ref{tab:cd_gmres} shows execution times for various values of the restart parameter of GMRES.
\begin{table}[htbp]
\caption{Varying the restart parameter of GMRES for the 3D convection-diffusion test case \eqref{eq:cd}, problem size $n = 100^{3}$.}
{
\begin{tabular}{ccccccc}
\hline
  $p$ & \multicolumn{3}{c}{$48$} & \multicolumn{3}{c}{$192$}\\
  \hline
  Restart & Clock (sec) & $k$ & $r$ & Clock (sec) & $k$ & $r$\\
  \hline
  5 & 344 & 917 & 9.98E-07 & 187 & 917 & 9.98E-07\\
  10 & 251 & 489 & 9.70E-07 & 149 & 489 & 9.70E-07\\
  20 & 274 & 318 & 9.44E-07 & 161 & 318 & 9.44E-07\\
  30 & 427 & 349 & 9.77E-07 & 247 & 349 & 9.77E-07\\
  40 & 614 & 385 & 9.65E-07 & 349 & 385 & 9.65E-07\\
  50 & 748 & 393 & 9.59E-07 & 440 & 393 & 9.59E-07\\
  100 & 1765 & 457 & 9.80E-07 & 969 & 457 & 9.80E-07\\
  (Full) & 2695 & 281 & 8.56E-07 & 1677 & 281 & 8.56E-07\\
\hline
\end{tabular}
}
\label{tab:cd_gmres}
\end{table}
This allows us to choose the value 10 as the experimentally optimal one, however, performances for a restart value of 20 were quite similar.

We therefore looked for performance variation of HSS(CG, GMRES(10)) according to its parameter $\alpha$ and the inner residual threshold $\varepsilon_{\text{in}}$ set for both CG and GMRES(10). Convergence was obtained from $\varepsilon_{\text{in}} = 10^{-2}$, which also demonstrated more efficiency than lower thresholds, as shown in Table \ref{tab:cd_ihss}.
\begin{table}[htbp]
\caption{Varying the parameter $\alpha$ and the inner residual threshold $\varepsilon_{\text{in}}$ of HSS(CG,GMRES(10)) for the 3D convection-diffusion test case \eqref{eq:cd}, problem size $n = 100^{3}$, number of processes $p = 192$.}
{
\begin{tabular}{cccccccccc}
\hline
  \multicolumn{5}{c}{$\varepsilon_{\text{in}}$ = 1.00E-02} & \multicolumn{5}{c}{$\varepsilon_{\text{in}}$ = 1.00E-06}\\
  \hline
  $\alpha$ & Clock (sec) & $k$ & $k_{\text{in}}$ & $r$ & $\alpha$ & Clock (sec) & $k$ & $k_{\text{in}}$ & $r$\\
  \hline
  0.7 & 718 & 213 & 2182 & 9.84E-07 & 0.9 & 2431 & 270 & 7331 & 9.85E-07\\
  0.6 & 712 & 186 & 2124 & 9.57E-07 & 0.8 & 2395 & 240 & 7129 & 9.85E-07\\
  0.5 & 665 & 162 & 1949 & 9.94E-07 & 0.7 & 2398 & 210 & 6986 & 9.84E-07\\
  0.4 & 844 & 164 & 2148 & 9.76E-07 & 0.6 & 2450 & 180 & 6916 & 9.84E-07\\
\hline
\end{tabular}
}
\label{tab:cd_ihss}
\end{table}
Quite surprisingly, the number of outer iterations even slightly increased when switching from $10^{-2}$ to $10^{-6}$.

While a restart value of 10 resulted in the most efficient executions of the GMRES solver, it does not necessarily prove to be the best choice for HSS(CG, GMRES(restart)) as well. Handling a combination of three parameters, $\alpha$, $\varepsilon_{\text{in}}$ and GMRES' restart, is clearly a major drawback of HSS(CG, GMRES(restart)), especially if, additionally, the number of processes (and so, possibly, the load per process) might have an impact too. Our two-stage-splitting-based HSS(${M}^{-1}$, ${F}^{-1}$) with single inner iteration takes the set of parameters back to $\alpha$, as in the case of exact HSS. Moreover, as mentioned in Section \ref{sec:impl}, avoiding inner solver function calls and loops might constitute an attractive feature, considering pure computer science aspects. This is shown here by comparing Tables \ref{tab:cd_ihss} and \ref{tab:cd_rihss}.
\begin{table}[htbp]
\caption{Varying the parameter $\alpha$ of HSS(${M}^{-1}$, ${F}^{-1}$) for the 3D convection-diffusion test case \eqref{eq:cd}, problem size $n = 100^{3}$.}
{
\begin{tabular}{ccccccc}
\hline
  $p$ & \multicolumn{3}{c}{$48$} & \multicolumn{3}{c}{$192$}\\
  \hline
  $\alpha$ & Clock (sec) & $k$ & $r$ & Clock (sec) & $k$ & $r$\\
  \hline
  6.0 & 566 & 2348 & 9.98E-07 & 252 & 2307 & 9.98E-07\\
  5.0 & 485 & 2008 & 9.99E-07 & 214 & 1965 & 9.94E-07\\
  4.0 & 399 & 1657 & 9.94E-07 & 177 & 1611 & 9.98E-07\\
  3.0 & 311 & 1288 & 9.90E-07 & 136 & 1239 & 9.70E-07\\
\hline
\end{tabular}
}
\label{tab:cd_rihss}
\end{table}
For $p = 192$ processes, best execution times of HSS(CG, GMRES(10)) and HSS(${M}^{-1}$, ${F}^{-1}$) are, respectively, 665 and 136 seconds. Note that the former performed 1949 inner iterations while the latter converged in 2576 inner iterations (2 $\times$ 1288 outer iterations since there is one inner iteration using ${M}^{-1}$ and another one using ${F}^{-1}$). Such a surprisingly quite small gap in convergence speed confirms the possibility to achieve a faster solver in execution time by avoiding inner function calls and loops. Still, an important drawback for HSS(${M}^{-1}$, ${F}^{-1}$) is that it turned divergent for $\alpha \le 2.0$.

Finally, Table \ref{tab:cd_aihss} shows that $\alpha = 3.0$ was experimentally optimal for the asynchronous HSS(${M}^{-1}$, ${F}^{-1}$) too. And here as well, divergence has been observed for $\alpha \le 2.0$.
\begin{table}[htbp]
\caption{Varying the parameter $\alpha$ of asynchronous HSS(${M}^{-1}$, ${F}^{-1}$) for the 3D convection-diffusion test case \eqref{eq:cd}, problem size $n = 100^{3}$.}
{
\begin{tabular}{ccccccccc}
\hline
  $p$ & \multicolumn{4}{c}{$48$} & \multicolumn{4}{c}{$192$}\\
  \hline
  $\alpha$ & Clock (sec) & $k_{\text{min}}$ & $k_{\text{max}}$ & $r$ & Clock (sec) & $k_{\text{min}}$ & $k_{\text{max}}$ & $r$\\
  \hline
  6.0 & 24 & 3134 & 4609 & 4.32E-07 & 7.46 & 7299 & 9491 & 4.83E-07\\
  5.0 & 22 & 2812 & 3969 & 4.31E-07 & 7.04 & 6832 & 9175 & 6.57E-07\\
  4.0 & 20 & 2573 & 3695 & 4.21E-07 & 6.82 & 6668 & 8846 & 5.12E-07\\
  3.0 & 17 & 2278 & 3080 & 5.49E-07 & 6.24 & 5950 & 7996 & 9.78E-07\\
\hline
\end{tabular}
}
\label{tab:cd_aihss}
\end{table}

\subsubsection{Performance comparison}

Using experimentally obtained optimal parameters, a performance comparison on $p = 48$ to $p = 192$ cores is summarized here in Table~\ref{tab:1}, where we dropped off the HSS(CG, GMRES(10)) due to memory limits exceeded for $p \le 120$.
\begin{table}
\caption{Performances from the 3D convection-diffusion test case \eqref{eq:cd}, problem size $n = 100^{3}$.}
{
\begin{tabular}{ccccccccccc}
  \hline
  & \multicolumn{3}{c}{GMRES(10)} & \multicolumn{3}{c}{HSS(${M}^{-1}$, ${F}^{-1}$, 3.0)} & \multicolumn{4}{c}{Async. HSS(${M}^{-1}$, ${F}^{-1}$, 3.0)} \\
  \hline
  $p$ & Clock & $k$ & $r$ & Clock & $k$ & $r$ & Clock & $k_{\text{min}}$ & $k_{\text{max}}$ & $r$ \\
  & (sec) & & & (sec) & & & (sec) & & & \\
  \hline
  48 & 251 & 489 & 9.70E-07 & 311 & 1288 & 9.90E-07 & 17 & 2278 & 3080 & 5.49E-07 \\
  72 & 197 & 489 & 9.70E-07 & 222 & 1222 & 9.92E-07 & 12 & 3401 & 3912 & 8.44E-07 \\
  96 & 239 & 489 & 9.70E-07 & 203 & 1177 & 9.92E-07 & 14 & 5682 & 6678 & 9.21E-07 \\
  120 & 151 & 489 & 9.70E-07 & 193 & 1228 & 9.97E-07 & 12 & 6541 & 8233 & 8.79E-07 \\
  144 & 169 & 489 & 9.70E-07 & 179 & 1229 & 9.93E-07 & 10 & 7176 & 9394 & 9.50E-07 \\
  168 & 150 & 489 & 9.70E-07 & 133 & 1240 & 9.89E-07 & 6.20 & 5526 & 7562 & 8.59E-07 \\
  192 & 149 & 489 & 9.70E-07 & 136 & 1239 & 9.70E-07 & 6.24 & 5950 & 7996 & 9.78E-07 \\
  \hline
\end{tabular}
}
\label{tab:1}
\end{table}
One can see a significant gain by asynchronous HSS(${M}^{-1}$, ${F}^{-1}$, 3.0), which was, e.g., at $p = 192$ processor cores, about 20 times faster (in execution time) than both GMRES(10) and synchronous HSS(${M}^{-1}$, ${F}^{-1}$, 3.0). While the second-stage splittings using preconditioners ${M}^{-1}$ and ${F}^{-1}$ were introduced here to achieve a fully asynchronous version of HSS, such a gap between the performances of synchronous and asynchronous HSS(${M}^{-1}$, ${F}^{-1}$, 3.0) in a homogeneous high-speed computational environment shows that there is a true advantage in resorting to asynchronous iterations, which is not due to possible programming biases introduced by this particular implementation of HSS.

\subsection{Results on the 2D structural dynamics problem}

\subsubsection{Optimal parameters}

The complex 2D structural dynamics test case \eqref{eq:sd} was run on an obtained discrete problem with $n = 350^{2}$ unknowns, using from $p = 24$ to $p = 54$ processor cores (one MPI process per core).

Table \ref{tab:sd_gmres} shows execution times for various values of the restart parameter of GMRES.
\begin{table}[htbp]
\caption{Varying the restart parameter of GMRES for the 2D structural dynamics test case \eqref{eq:sd}, problem size $n = 350^{2}$, number of processes $p = 48$.}
{
\begin{tabular}{cccc}
\hline
  Restart & Clock (sec) & $k$ & $r$\\
  \hline
  5 & 5405 & 36594 & 1.00E-06\\
  10 & 3960 & 19679 & 1.00E-06\\
  20 & 3068 & 9072 & 1.01E-06\\
  30 & 3053 & 6386 & 1.02E-06\\
  40 & 3158 & 5125 & 1.04E-06\\
  50 & 3084 & 4080 & 9.84E-07\\
  100 & 3433 & 2727 & 7.89E-07\\
  (Full) & 7898 & 789 & 9.63E-07\\
\hline
\end{tabular}
}
\label{tab:sd_gmres}
\end{table}
This allows us to choose the value 30 as the experimentally optimal one, however, performances for restart values of 20 to 50 were quite similar.

Both HSS(CG, GMRES(30)) and HSS(${M}^{-1}$, ${F}^{-1}$) failed to converge within two hours of execution on $p = 48$ cores for various values of their parameters, which made them unpractical for the current test case.

Nevertheless, asynchronous HSS(${M}^{-1}$, ${F}^{-1}$) took reasonable times to converge, and Table \ref{tab:sd_aihss} shows an experimentally optimal $\alpha = 2.0$. Divergence was observed for $\alpha \le 1.0$.
\begin{table}[htbp]
\caption{Varying the parameter $\alpha$ of asynchronous HSS(${M}^{-1}$, ${F}^{-1}$) for the 2D structural dynamics test case \eqref{eq:sd}, problem size $n = 350^{2}$, number of processes $p = 48$.}
{
\begin{tabular}{ccccc}
\hline
  $\alpha$ & Clock (sec) & $k_{\text{min}}$ & $k_{\text{max}}$ & $r$\\
  \hline
  5.0 & 273 & 398754 & 493820 & 7.19E-07\\
  4.0 & 235 & 349111 & 425328 & 8.71E-07\\
  3.0 & 198 & 293439 & 357005 & 1.04E-06\\
  2.0 & 156 & 231787 & 281838 & 9.50E-07\\
\hline
\end{tabular}
}
\label{tab:sd_aihss}
\end{table}

\subsubsection{Performance comparison}

Using experimentally obtained optimal parameters, a performance comparison on $p = 24$ to $p = 54$ cores is summarized in Table~\ref{tab:2}.
\begin{table}
\caption{Performances from the complex 2D structural dynamics test case \eqref{eq:sd}, problem size $n = 350^{2}$.}
{
\begin{tabular}{cccccccc}
  \hline
  & \multicolumn{3}{c}{GMRES(30)} & \multicolumn{4}{c}{Async. HSS(${M}^{-1}$, ${F}^{-1}$, 2.0)}\\
  \hline
  $p$ & Clock (sec) & $k$ & $r$ & Clock (sec) & $k_{\text{min}}$ & $k_{\text{max}}$ & $r$ \\
  \hline
  24 & 2941 & 6486 & 9.99E-07 & 308 & 183861 & 203002 & 8.50E-07 \\
  30 & 2722 & 6419 & 9.99E-07 & 253 & 212597 & 249716 & 8.81E-07 \\
  36 & 2967 & 6510 & 1.02E-06 & 241 & 236977 & 277301 & 9.86E-07 \\
  42 & 2656 & 6479 & 1.02E-06 & 154 & 211052 & 257389 & 1.01E-06 \\
  48 & 3053 & 6386 & 1.02E-06 & 156 & 231787 & 281838 & 9.50E-07 \\
  54 & 2829 & 6479 & 1.01E-06 & 159 & 251221 & 310456 & 9.13E-07 \\
  \hline
\end{tabular}
}
\label{tab:2}
\end{table}
Again, a significant gain is obtained by asynchronous HSS(${M}^{-1}$, ${F}^{-1}$, 2.0), which was, e.g., at $p = 48$ processor cores, about 20 times faster than GMRES(30), similarly to the real 3D convection-diffusion test case. Here as well an even more important performance gap is observed between asynchronous and synchronous HSS(${M}^{-1}$, ${F}^{-1}$, 2.0) which did not terminate within 7200 seconds. This confirms, for the complex test case as well, the benefit purely from asynchronous iterations. 

}

\section{Conclusion}\label{sec:con}

Asynchronous alternating iterations are revealed here as a practical breakthrough in improving computational time of parallel solution of non-Hermitian problems, compared to the well-known GMRES and HSS methods. Classical asynchronous convergence conditions are investigated for a general practical parallel scheme of alternating iterations. In particular, it can result in a two-stage variant of the HSS method with one inner iteration for each of the outer alternating ones. Performance experiments have been conducted for such an asynchronous variant which has significantly outperformed both the GMRES and the classical HSS methods, both on a real convection-diffusion and a complex structural dynamics problem.

\section*{Acknowledgement}

The paper has been prepared with the support of the ``RUDN University Program 5-100'', the French national program LEFE/INSU, the project ADOM (M\'ethodes de d\'ecomposition de domaine asynchrones) of the French National Research Agency (ANR), and using HPC resources from the ``M\'esocentre'' computing center of CentraleSup\'elec and \'Ecole Normale Sup\'erieure Paris-Saclay supported by CNRS and R\'egion \^Ile-de-France.

\bibliography{ref}
\bibliographystyle{abbrv}

\end{document}